\documentclass[12pt]{amsproc}

\def\ThmNum{yes}

\newfont{\bbf}{msbm10 at 12pt}
\newfont{\bbfsm}{msbm10 at 9pt}

\def\N{\mathbb {N}}

\def\Z{\mathbb {Z}}
\def\R{\mathbb {R}}
\def\C{\mathbb {C}}
\def\Cbar{\overline{\C}}

\def\phi1{\phi}
\def\phi{\varphi}
\def\eps{\varepsilon}
\def\theta{\vartheta}
\def\Re{\mbox{\rm Re}}

\def\sm{\setminus}

\ifx\ThmNum\undefined
      \newtheorem{theorem}{Theorem}[section]
      
\else
      \newtheorem{theorem}{Theorem}
      
\fi

\newtheorem{lemma}[theorem]{Lemma}
\newtheorem{definition}[theorem]{Definition}

\newtheorem{corollary}[theorem]{Corollary}

\newtheorem{example}[theorem]{Example}

\def\proof{\par\medskip\noindent {\sc Proof. }}
\def\proofof #1 {\par\medskip\noindent {\sc Proof of #1. }}

\def\sketchof #1 {\par\medskip\noindent {\sc Sketch of proof of #1. }}
\def\Box{\framebox[10pt]{\rule{0pt}{3pt}}}
\def\nix{\rule{0pt}{2pt}}
\def\qed{\qedd\par\medskip\noindent}
\def\qedd{\nix\nolinebreak\hfill\hfill\nolinebreak$\Box$}
\def\remark{\par\medskip \noindent {\sc Remark. }\nopagebreak}
\def\lineclear
       {\rule{0pt}{0pt}\nopagebreak\par\nopagebreak\noindent}

\def\reminder #1 {{\sf #1}}
\def\hide #1 {}

\newcommand{\sumr}{\mathop{\to\!\!\!\!\!\!\!\!\sum}\limits}
\newcommand{\suml}{\mathop{\leftarrow\!\!\!\!\!\!\!\!\sum}\limits}
\newcommand{\sumrl}{\mathop{\leftrightarrow\!\!\!\!\!\!\!\!\sum}\limits}
\newcommand{\prodr}{\mathop{\to\!\!\!\!\!\!\!\prod}\limits}

\def\up#1#2{
\rule{0pt}{9pt}^{#2}\!#1
}
\def\strut{\rule{0pt}{11pt}}

\def\remarks{\par\medskip \noindent {\sc Remarks. }\nopagebreak}

\title{Fractional Sums and Euler-like Identities}

\author{Markus M\"uller}
\address{
Institut f\"ur Mathematik, Technische Universit\"at Berlin,
Stra{\ss}e des 17. Juni 135, D-10623 Berlin, Germany
}
\email{mueller@math.tu-berlin.de}
\author{Dierk Schleicher}
\address{
School of Engineering and Science, IUB: International
University Bremen, Postfach 750~561, D-28725~Bremen, Germany}
\email{dierk@iu-bremen.de}

\keywords{Fractional Sum, summation, interpolation, summation identities}
\subjclass[2000]{33-99, 40-99, 40A25, 41A05}

\begin{document}

\renewcommand{\thefootnote}{}
\footnotetext{{\em Date:} October 16, 2006}
\renewcommand{\thefootnote}{\arabic{footnote}}

\begin{abstract}
We introduce a natural definition for sums of the form
\[
\sum_{\nu=1}^x f(\nu)
\]
when the number of terms $x$ is a rather arbitrary real or even complex number.
The resulting theory includes the known interpolation of the
factorial by the $\Gamma$ function or Euler's little-known formula
\(
\sum_{\nu=1}^{-1/2} \frac 1\nu = -2\ln 2
\,\,.
\)

Many classical identities like the geometric
series and the binomial theorem nicely extend to this more general setting.
Sums with a fractional number of terms are closely related to special 
functions, in
particular the Riemann and Hurwitz $\zeta$ functions. A number of results about
fractional sums can be interpreted as classical infinite sums or products or as
limits, including identities like
\hide{
\[
   \lim_{n\to\infty}\left[
      e^{\frac n 8 (4n+1)} n^{-\frac 1 {16}-\frac 1 2 n(n+1)} (2 \pi)^{-\frac n 4}
      \prod_{k=1}^n \frac {\Gamma(k+1)^k}{\Gamma(k+\frac 1 2)^{(k-\frac 1 2)}}
   \right]
   =(2e)^{\frac 1 {24}} A^{\frac 3 4} \exp\left(-\frac{7\zeta(3)}{32\pi^2}\right)
\]
}
\[  \lim_{n\to\infty}\left[     e^{\frac n 4(4n+1)}n^{-\frac 1 8 - n(n+1)}(2\pi)^{-\frac n 2}     \prod_{k=1}^{2n} \Gamma\left(1+\frac k 2\right)^{k(-1)^k}  \right]\hide{  =(2e)^{\frac 1 {12}} A^{\frac 3 2} \exp\left(-\frac{7\zeta(3)}{16\pi^2}\right)}
= \sqrt[12]{2} \exp\left(\frac{5}{24} - \frac 3 2 \zeta'(-1) -\frac{7\zeta(3)}{16\pi^2}\right)\]
some of which seem to be new.

\end{abstract}

\maketitle


\bigskip
\begin{center}
\begin{minipage}{120mm}
\def\tocname{Contents}
\tableofcontents
\end{minipage}
\end{center}

\section{Introduction}
\label{SecIntro}

Sums of the form $\sum_{\nu=1}^x f(\nu)$ are defined
classically only when the number of terms $x$ is a positive
integer or $\infty$.
There have been certain attempts to interpolate summations to
non-integer numbers of terms; the most famous one is probably the
interpolation of the factorial, which after taking logarithms can
be written as
\[
\sum_{\nu=1}^x \ln\nu = \ln\Gamma(x+1)
\]
using the well-known $\Gamma$ function which was introduced
for this very purpose (compare Example~\ref{ExFactorial}).

We propose a systematic way to extend summations to non-integer
numbers of terms:
there is a natural and essentially unique way to do this starting
from the continued summation property
$\sum_{a}^b+\sum_{b+1}^c=\sum_a^c$.
The works of Euler~\cite[p.~97]{Euler} and
Ramanujan~\cite[Chapter~6, entry 4(i)]{B}
show that both of them had also looked into a similar direction:
for example, Euler~\cite[p.~101]{Euler} has the formula
$\sum_{\nu=1}^{-1/2}\nu^{-1}=-2\ln 2$ (see
Corollary~\ref{CorRiemannHurwitz2}). However, we are not aware of any
attempts of a treatment beyond sporadic examples.

In this paper, we define ``fractional sums'', explore the
consequences of this definition and derive several of its
properties. While we give a precise definition in Section~\ref{SecDef},
a motivating special case follows from continued summation with
summation boundary at infinity: the identity
\[
\sum_{\nu=1}^x +
\sum_{\nu=x+1}^\infty
=
\sum_{\nu=1}^\infty
\]
certainly holds for $x\in\N$. The sum on the right makes sense 
classically, and so
does the middle one: $\sum_{\nu=x+1}^\infty 
f(\nu)=\sum_{\nu=1}^\infty f(\nu+x)$.
We can thus use this equation to define the left sum for $x\in\C$. Our
Definition~\ref{DefFracSummable} extends this simple idea to a larger class of
functions. The general idea is to shove the ill-defined terms to 
$+\infty$ where
they can be estimated precisely in the limit. In 
Section~\ref{SecLeft}, we explore
the related definition when the limit $-\infty$ rather than $+\infty$ is used.

We generalize some well-known algebraic identities from an integer 
number of terms
to a fractional (even complex) number of terms
(Sections~\ref{SecCharacterization}, \ref{SecProducts} and \ref{SecLeft}). For
example, the finite geometric series (for $x\in\N$)
\[
\sum_{\nu=0}^x q^\nu = \frac{q^{x+1}-1}{q-1}
\]
nicely generalizes to real and even complex values of $x$
(Theorem~\ref{ThmGeometricSeries}). Similarly, there is the
well-known formula for squaring a finite sum: $(a+b)^2=a^2+b^2+2ab$, or more
generally for $x\in\N$:
\begin{equation}
\left(\sum_{\nu=1}^x a_\nu\right)^2
=
\sum_{\nu=1}^x a^2_\nu + 2\sum_{\mu<\nu}a_\mu a_\nu \,\,;
\label{EqSquaring}
\end{equation}
in Corollary~\ref{CorQuadrature2}, we prove an analogous formula for 
arbitrary $x\in
\C$.

While these identities express properties of our ``fractional sums'', they have
special cases which can be rewritten in terms of classical infinite sums.
This way, we obtain several identities involving infinite sums:
some of them are well known, while others seem to be new, and many others can
be produced in a similar spirit. We give a few examples in
Section~\ref{SecEvalInfinite}, including a proof of the formula
\[  \lim_{n\to\infty}\left[     e^{\frac n 4(4n+1)}n^{-\frac 1 8 - n(n+1)}(2\pi)^{-\frac n 2}     \prod_{k=1}^{2n} \Gamma\left(1+\frac k 2\right)^{k(-1)^k}  \right]\hide{  =(2e)^{\frac 1 {12}} A^{\frac 3 2} \exp\left(-\frac{7\zeta(3)}{16\pi^2}\right)}
= \sqrt[12]{2} \exp\left(\frac{5}{24} - \frac 3 2 \zeta'(-1) -\frac{7\zeta(3)}{16\pi^2}\right)\]
which we have not found anywhere else. 
\hide{
($A=\exp\left(-\frac 1 {12}-\zeta'(-1)\right)$ denotes the Glaisher-Kinkelin constant.)
}

The occurrence of special functions like $\Gamma$ or Riemann's $\zeta$ function
in this example is no coincidence; our definition of
fractional sums leads in many cases quite naturally to special 
functions. Some of
the relations to $\zeta$ functions are discussed in 
Section~\ref{SecRiemannHurwitz}.

Some results in this paper have been announced in \cite{FracSumsProc}.

We should mention that the well-known Euler-Maclaurin formula can 
also be used to
estimate sums with non-integer numbers of terms provided the number 
$x$ of terms is
large. In many cases, the difference between the values of 
Euler-Maclaurin and our
method tends to zero as $x\to\infty$ through the reals. However, the 
methods are
fundamentally different: for example, our method does not require the 
integrand to
be differentiable or even continuous.

This paper grew out of a project at the German youth science fair
``Jugend forscht'' by the first author when he was a high school
student (unaware of the works of Euler and Ramanujan).
Both of us are no experts in this area, and we apologize if references or
due credit are missing.

{\bf Acknowledgments.} We have received encouragement,
suggestions and support from a number of people, including
Richard Askey, Mourad Ismail, Irwin Kra, John Milnor, Michael Stoll and the
late Judita Cofman, as well as the audiences in M\"unchen, Stony
Brook and Bremen, and in Co\-pen\-hagen on the OPSFA conference.
We would like to thank them all.

\section{The Fundamental Definition}
\label{SecDef}

We begin by a natural definition for polynomials.
\begin{definition}[Fractional Sums for Polynomials]
\label{DefFracSumsPoly}\lineclear
For a polynomial $p:\C\to\C$, let $P$ be the unique polynomial with
$\sum_{\nu=1}^n p(\nu)=P(n)$ for all $n\in\N$. Then we define for 
every $x\in\C$
\[
      \sumrl_{\nu=1}^x p(\nu):=P(x)\,\,.
\]
Moreover, for arbitrary $a,b\in\C$, we define
\[
      \sumrl_{\nu=a}^b p(\nu):=\sumrl_{\nu=1}^{b-a+1} p(\nu+a-1)
=P(b)-P(a-1)
\,\,.
\]
\end{definition}

In this paper, we extend this concept to a larger class of functions 
as follows.

\begin{definition}[Approximate Polynomial]
\label{DefApproxPoly} \lineclear
Let $U\subset \C$ and $\sigma\in\N\cup\{-\infty\}$.
A function $f\colon U\to\C$ will be called a {\em (right) approximate 
polynomial of
degree $\sigma$} if the following conditions are satisfied:
\begin{itemize}
\item
all $u\in U$ satisfy $u+1\in U$ \,\,;
\item
there exists a sequence of polynomials $(p_n)_{n\in\N}$ of fixed
degree $\sigma$ such that for every $x\in U$,
\[
\left| f(n+x) - p_n(n+x) \right| \longrightarrow 0
\quad \mbox{as $n\to+\infty$} \,\,.
\]
\end{itemize}
\end{definition}

This is a semi-local condition and not too restrictive; only the 
behavior of $f(x)$
as $\Re(x)\to +\infty$ matters. For example, every $f\colon\C\to\C$ 
with $f(x)\to
0$ as $\Re(x)\to+\infty$ is approximately polynomial of degree 
$-\infty$, and the
functions $f(x)=\ln x$ and $f(x)=\sqrt x$ on $\R^+$ are approximately 
polynomial
of degree $0$ (i.e.\ approximately constant). The class of 
approximate polynomials
is large enough for many interesting applications.

Now comes our general definition of fractional sums. It uses the
approximating polynomials as well as their fractional sums as
defined above. After the formal definition, we try to motivate this 
definition and
explain it in a number of special cases.

\begin{definition}[Fractional Sum and Product]
\label{DefFracSummable} \lineclear
An approximate polynomial $f\colon U\to \C$ of degree 
$\sigma\in\N\cup\{-\infty\}$ will be called
{\em right summable} \hide{ of degree $\sigma$}
if for every $a,b+1\in U$, the limit
\[
     \lim_{n\to\infty}\left(
        \sumrl_{\nu=n+a}^{n+b}p_n(\nu)+\sum_{\nu=1}^n \left(\strut
           f(\nu+a-1)-f(\nu+b)
        \right)
     \right)
\]
exists.
In this case, this limit will be the definition for the fractional 
sum of $f$ from $a$ to $b$; we denote it by
\[
\sumr_{\nu=a}^b f(\nu)
\qquad\mbox{or briefly}\qquad
\sumr_a^b f
\,\,.
\]
Moreover, we can define fractional products by
\[
     \prodr_{\nu=a}^b f(\nu):=\exp\left(
        \sumr_{\nu=a}^b \ln f(\nu)
     \right)\,\,,
\]
whenever $\ln f$ is right summable.
\end{definition}

\remarks \nopagebreak
\begin{itemize}
\item
In the limit, $n\in\N$ is always taken to be an integer.
\item
The value of the sum is independent of the choice of the approximating
polynomials $p_n$.
\item
If $b-a\in\N$, then the limit exists and agrees with the classical value
of the sum.

\item
If $f$ is a polynomial, then Definition~\ref{DefFracSummable} is 
consistent with
Definition~\ref{DefFracSumsPoly}, hence $\sumr f = \sumrl f$ for 
arbitrary complex
summation boundaries.

\item
We use the notation $\sumr$ for our fractional sum using ``right summable
functions'' (where the undefined terms are shoved to $+\infty$); similarly, in
Section~\ref{SecLeft}, we introduce the symbol $\suml$ for ``left summable
functions''. The symbol $\sumrl$ is used for polynomials where both concepts
coincide trivially.

\item
It may be helpful to write down the first few cases of $\sigma$
explicitly (for $a=1$ and $b=x$); the cases $\sigma\geq 1$ involve
some calculations.

\begin{eqnarray*}
\sigma=-\infty: &&
\sumr_{\nu=1}^xf(\nu)
=\sum_{\nu=1}^{\infty}\left(f(\nu)-f(\nu+x)\strut\right)
\label{EqSigma-Inf}
\\
\sigma=0: &&
\sumr_{\nu=1}^x f(\nu) =
\lim_{n\to+\infty}\left(
x\cdot f(n)+\sum_{\nu=1}^n \left(f(\nu)-f(\nu+x) \strut\right)
\right)
\label{EqSigmaZero}
\\
\sigma=1: &&
\sumr_{\nu=1}^x f(\nu) =
\lim_{n\to+\infty}\left(
x\cdot f(n)+\frac{x(x+1)}{2} \cdot\left(f(n+1)-f(n) \strut\right)
\rule{0pt}{18pt}\right. \nonumber
\\
&& \qquad\qquad\qquad\qquad\qquad \left.
+\sum_{\nu=1}^n \left(f(\nu)-f(\nu+x) \strut\right) \right)
\,\,.
\label{EqSigmaOne}
\end{eqnarray*}

\end{itemize}

We will now try to motivate our definition. We begin with the case
$\sigma=-\infty$, so that all $p_n\equiv 0$. In this case, we calculate the sum
of $f$ from $1$ to $x+n$ in two ways
\[
\sumr_{1}^{x+n}f =
\sumr_{1}^x f + \sum_{x+1}^{x+n} f = \sum_{1}^n f + \sumr_{n+1}^{x+n} f \,\,,
\]
so we obtain
\begin{equation}
\label{EqHeuristics2}
\sumr_{1}^x f = \sum_{1}^n f + \sumr_{n+1}^{x+n} f - \sum_{x+1}^{x+n} f
\,\,.
\end{equation}
If  $f(x)$ tends to $0$ as $x\to\infty$, then it makes sense to require that
$\sumr_{n+1}^{x+n}\to 0$ as $n\to\infty$; the remaining two terms on 
the right hand
side exactly yield the definition above. This motivates the definition for
$\sigma=-\infty$. Note that in the last sum, the difference between 
upper and lower
boundaries is an integer, so this sum is defined classically. In this 
heuristics,
sums with a non-integer number of terms are denoted by our generalized symbol
$\sumr$.

The next case is that of ``approximately constant'' functions $f$: suppose that
$f=\ln$. In this case, the function $f$ can be approximated by constants over
regions of bounded diameters: for every bounded domain $W$ and every
$\eps>0$, there is an $N\in\N$ such that for all $n\geq N$ there is a 
$C_n\in\C$
for which $|f(z+n)-C_n|<\eps$ uniformly for $z\in W$. It thus
makes sense to estimate $\sumr_{\nu=n+1}^{n+x} f(\nu)\approx x
f(n)$ (the approximating constants change with $n$, but the
quality of the approximation improves for large $n$). This leads to the case
$\sigma=0$ of our definition. Note that exponentiating this example immediately
leads to the interpolation of the factorial by the $\Gamma$ function; 
see below.

The general case is similar: in the right hand side of 
(\ref{EqHeuristics2}), the
first and last sums always lead to the same term $\sum_{\nu=1}^n
(f(\nu+1)-f(\nu+x))$, while the second sum can be estimated well by the exact
formulas for polynomials.

Note that ``approximately polynomial of degree $-\infty$'' is a more general
condition than $f(x)\to 0$ as $x\to\infty$; similarly, 
``approximately constant''
is more general than the uniform condition given above for $f=\ln$:
the ``approximately polynomial'' condition evaluates $f$ only at $\N$ 
and at $x+\N$,
which are the only values used in our definition.

We can now show that our definition fits in with the known interpolation of the
factorial function.
\begin{example}[The Extended Factorial]
\label{ExFactorial} \lineclear
For every $x\in\C\sm\{-1,-2,-3\dots\}$, the factorial has the following
product formula
\[
\prodr_{\nu=1}^x \nu =
\lim_{n\to+\infty}\left( n^x \prod_{\nu=1}^n\frac{\nu}{\nu+x}\right)
=\Gamma(x+1) \,\,.
\]
\end{example}
\proof
As noted above, $\ln$ is approximately constant ($\sigma=0$), so our
definition reads
\begin{equation}
\sumr_{\nu=1}^x \ln \nu = \lim_{n\to\infty}\left(
x\ln n
+ \sum_{\nu=1}^n \left(\ln\nu-\ln(\nu+x) \strut\right)
\right)
\label{EqDefLogSum}
\end{equation}
up to an additive term in $2\pi i\Z$. For every finite $n$, the sum is
well-defined provided $x\in\C\sm\{-1,-2,-3\dots\}$. For $\nu$ large,
$\ln(\nu)-\ln(\nu+x)=\ln(\nu/(\nu+x))$ is evaluated as the principal
branch, and finitely many choices of the branch for small values of
$\nu$ are irrelevant for convergence. For $\ln n$, we use the principal
branch.

It is readily verified that the limit in (\ref{EqDefLogSum}) exists, 
so the sum is
well defined (up to finitely many additive summands $2\pi i$, which 
are canceled by
the subsequent exponentiation). By definition of the product in
Definition~\ref{DefFracSummable}, we get
\begin{eqnarray*}
\prodr_{\nu=1}^x \nu
&=&
\exp\left(\sumr_{\nu=1}^x \ln \nu\right)
= \lim_{n\to\infty}\exp\left(
x\ln n + \sum_{\nu=1}^n \left(\strut\ln\nu-\ln(\nu+x)\right) \right) \\
&=&
\lim_{n\to\infty}\left(
n^x \prod_{\nu=1}^n \frac{\nu}{\nu+x}\right) \,\,.
\end{eqnarray*}
It is well known \cite[Gauss' Formula VII.7.6]{Conway} that this equals
$\Gamma(x+1)$.
\qed

The following basic properties follow immediately from the
definition:
\begin{theorem}[Basic Properties of Fractional Sums]
\label{ThmLinearity} \lineclear
Fractional sums have the following properties for arbitrary
$a,b,c,d\in\C$:
\begin{itemize}
\item
Linearity: \,\,\,$\displaystyle c\sumr_{\nu=a}^b f(\nu)+d\sumr_{\nu=a}^b g(\nu)
=\sumr_{\nu=a}^b\left(\strut
      c f(\nu)+d g(\nu)
\right)\,\,,$
\item
Continued Summation:\,\,\,$\displaystyle
\sumr_{\nu=a}^b f(\nu)+\sumr_{\nu=b+1}^c f(\nu)=\sumr_{\nu=a}^c f(\nu)\,\,,$
\item
Index Shifting:\,\,\,$\displaystyle \sumr_{\nu=a}^b
f(\nu+c)=\sumr_{\nu=a+c}^{b+c} f(\nu)\,\,,$
\end{itemize}
whenever two of the three fractional sums (in the last case: one of
the two sums) exist.
\qedd
\end{theorem}

Note also that  $\sumr_x^{x-1} f = 0$ and $\sumr_x^x f = f(x)$,
and more generally for $n\in\N$
\[
\sumr_{x}^{x+n} f = f(x)+f(x+1)+\ldots+f(x+n)\,\,.
\]

\section{Characterization of Fractional Sums}
\label{SecCharacterization}

The following result can be thought of as an analog to the 
fundamental theorem of
calculus. We use the notation $U^\pm:=\{u\pm 1\colon u\in U\}$ and $\Delta
S(x):=S(x)-S(x-1)$.

\begin{lemma}[Summation Formula For Approximate Polynomials]
\label{LemSummationApprox} \lineclear
Let $S\colon U\to\C$ be an approximate polynomial of degree
$\sigma\in\N\cup\{-\infty\}$, such that $0\in U$ (and thus $\N\subset U$ by definition).
Then $\Delta S\colon U^+\to\C$ is an approximate polynomial of degree 
$\sigma-1$;
moreover, for all $x\in U$, the sum $\displaystyle\sumr_{\nu=1}^x 
\Delta S(\nu)$
exists in $U$ and equals $S(x)-S(0)$.

Conversely, if $f\colon U^+\to \C$ is an approximate polynomial of 
degree $\sigma$
which is right summable, then the function $S(x):=\displaystyle\sumr_{\nu=1}^x
f(\nu)\colon U\to\C$ is approximately polynomial of degree $\sigma+1$ 
and satisfies
$\Delta S=f$.
\end{lemma}

\remark
The arithmetic of subtracting or adding degrees is like with 
differentiation and
integration of polynomials; in particular, speaking informally, we have
$0-1=-\infty$ and $-\infty-1=-\infty$, while $-\infty+1=0$ (the special case
$-\infty+1=-\infty$ is possible).

\proof
Set $f:=\Delta S$.
Let $P_n$ be the approximating polynomials of degree $\sigma$ for $S$ as in
Definition~\ref{DefApproxPoly}, and set $p_n(\nu):=P_n(\nu)-P_n(\nu-1)$.
For $x\in U^+$ and $n\in N$, we have
\begin{eqnarray*}
f(n+x)&=&
S(n+x)-S(n+x-1) \\
&=& P_n(n+x)-P_n(n+x-1)+o(1) \\
&=& p_n(n+x) +o(1) \,\,
\end{eqnarray*}
as $n\to\infty$,
and $f$ is approximately polynomial of degree $\sigma-1$.

For all $x\in U$ and $n\in\N$, we obviously have
\begin{equation}
S(n+x)
=
S(x)+\sum_{\nu=x+1}^{x+n} f(\nu)
\label{EqSummation1}
\end{equation}
and hence
\begin{eqnarray}
S(x)-S(0)
&=&
S(n+x)-S(n)-\sum_{\nu=x+1}^{x+n}f(\nu)+\sum_{\nu=1}^n f(\nu)
\label{EqSummation2}
\\
S(x)-S(0)&=& P_n(n+x)-P_n(n)+ o(1)
            +\sum_{\nu=1}^n \left(f(\nu)-f(\nu+x)\strut\right)
\,\,.
\label{EqSummation3}
\end{eqnarray}
We have $P_n(n+x)=P_n(n)+\sumrl_{\nu=n+1}^{n+x} p_n(\nu)$ for all
$x\in\C$, so we get
\begin{equation}
S(x)-S(0)=\lim_{n\to\infty}\left[\sumrl_{\nu=n+1}^{n+x} p_n(\nu) +
     \sum_{\nu=1}^n \left(f(\nu)-f(\nu+x)\strut\right) \right]
\label{EqSummation4}
\end{equation}
which is the definition of $\sumr_{\nu=1}^x f(\nu)$.

For the converse, we first observe that (\ref{EqSummation2}) holds 
trivially. Now
we read the proof backwards:
$S$ is defined via (\ref{EqSummation4}) because $f$ is right summable.
Let $p_n$ be the approximating polynomials of $f$ of degree $\sigma$. 
Then there
are polynomials $P_n$ of degree $\sigma+1$ such that
\[
\sumrl_{\nu=n+1}^{n+x} p_n(\nu) = P_n(n+x) - P_n(n)\,\,,
\]
which implies (\ref{EqSummation3}).
Combining (\ref{EqSummation2}) and (\ref{EqSummation3}), we get
\[
S(n+x)-S(n) = P_n(n+x)-P_n(n) + o(1) \,\,.
\]
Starting from $p_n$, the polynomials $P_n$ are defined only up to an
additive constant, which we may choose so that $P_n(n)=S(n)$. The claim
follows.
\qed

\section{Products of Fractional Sums}
\label{SecProducts}

In this section, we show that products like
$(a_1+a_2+\dots+a_n)\cdot(b_1+b_2+\dots+b_n)$ can be multiplied out for fractional
$n$ just like for integers. 

\begin{lemma}[Products of Fractional Sums]
\label{LemProdSums} \lineclear
Let $f,g\colon U\to\C$ be right summable functions such that
$x\mapsto \left(\sumr_{\nu=1}^x f(\nu)\right)\cdot\left(\sumr_{\nu=1}^x g(\nu)\right)$
is an approximate polynomial. Then every $x\in U$ satisfies
\begin{eqnarray*}
\left(\sumr_{\nu=1}^x f(\nu)\right)\cdot
\left(\sumr_{\nu=1}^x g(\nu)\right)
&=&
\sumr_{\nu=1}^x
\left(
f(\nu)g(\nu)
+
f(\nu)\cdot \sumr_{k=1}^{\nu-1}g(k)+
g(\nu)\cdot \sumr_{k=1}^{\nu-1}f(k)
\right)
\,\,.
\end{eqnarray*}
\end{lemma}
\proof
Let $F(x):=\sumr_{\nu=1}^x f(\nu)$ and $G(x):=\sumr_{\nu=1}^x g(\nu)$. Then $F(0)\cdot G(0)=0$ and
\begin{eqnarray*}
   \Delta\left(\strut F(x)G(x)\right)&\equiv &F(x)G(x)-F(x-1)G(x-1)\\
&=& \left(\strut f(x)+F(x-1)\right)\left(\strut g(x)+G(x-1)\right)
- F(x-1)G(x-1) \\
&=& f(x)g(x)+f(x)G(x-1)+g(x)F(x-1)\,\,.
\end{eqnarray*}
Since $F(x)G(x)$ is by assumption an approximate polynomial,
the first half of Lemma~\ref{LemSummationApprox} applies and proves the claim.
\qed

\begin{corollary}[Squares of Fractional Sums]
\label{CorQuadrature2} \lineclear
Suppose $f\colon U\to\C$ is a right summable function such that $x\mapsto\left(\sumr_{\nu=1}^x f(\nu)\right)^2$
is an approximate polynomial. Then
\[
\left[\sumr_{\nu=1}^x f(\nu)\right]^2 =
\sumr_{\nu=1}^x \left[
f^2(\nu)+2f(\nu)\sumr_{k=1}^{\nu-1}f(k)\right]
\,\,.
\]
\end{corollary}
\qed
\remark
Since $f$ is approximately polynomial, so is $F:=\sumr f$ by
Lemma~\ref{LemSummationApprox}; it is not automatic that $F^2$ is also
approximately polynomial, as the example $F(x)=x+\sin(x)/x$ shows.

\section{Left Summation and Binomial Series}
\label{SecLeft}

Classically, sums $\sum_{\nu=1}^N$ are defined only for integers $N$
with $N\geq 1$ or $N\geq 0$. If $N\in\N$, it will be natural for us to define
\begin{equation}
\sum_{\nu=1}^{-N} f(\nu):= -\sum_{\nu=-N+1}^0 f(\nu) \,\,;
\label{EqNegSum}
\end{equation}
this is the only way to extend the continued summation property from
Theorem~\ref{ThmLinearity}, and this works for {\em every} $f$ which is defined
on $\{-N+1,-N+2,\dots,0\}$. {\em Warning:} this is in contradiction to possible
conventions like
$\displaystyle \sum_{\nu=1}^{-N} = 0$ \quad or \quad
$\displaystyle\sum_{\nu=1}^{-N}=\sum_{\nu=-N}^1$.

Similarly, if $x-y\in-\N$, we set
\begin{equation}
\label{EqNegSum2}
\sum_{\nu=y}^x f(\nu):= -\sum_{\nu=x+1}^{y-1}f(\nu)
\end{equation}
for any $f$ which is defined at the finitely many points
$\{x+1,x+2,\dots,y-1\}$. With this convention,
Equation~(\ref{EqNegSum}) is valid for any integer $N$.

Our Definitions~\ref{DefApproxPoly} and \ref{DefFracSummable} of 
right approximate
polynomials and  $\sumr_{\nu=1}^x
f(\nu)$ use the behavior of $f$ as $\nu\to+\infty$; completely 
analogously, one can
also do this for $\nu\to-\infty$. The formal definition of {\em left
approximate polynomials} is analogous to Definition~\ref{DefApproxPoly}, except
that the limit $n\to+\infty$ is replaced by $n\to-\infty$.

\begin{definition}[Left Fractional Sum]
\label{DefFracSummableLeft} \lineclear
A left approximate polynomial $f\colon U\to \C$ of degree
$\sigma\in\N\cup\{-\infty\}$ will be called
{\em left summable}
if for every $a,b+1\in U$, the limit
\[
     \lim_{n\to-\infty}\left(
        \sumrl_{\nu=n+a}^{n+b}p_n(\nu)+\sum_{\nu=1}^n \left(\strut
           f(\nu+a-1)-f(\nu+b)
        \right)
     \right)
\]
exists.
In this case, this limit will be the definition for the left 
fractional sum of $f$
from $a$ to $b$; we denote it by
\[
\suml_{\nu=a}^b f(\nu)
\qquad\mbox{or briefly}\qquad
\suml_a^b f
\,\,.
\]
\end{definition}

\remark
In general, the existence of the two sums
$\suml_{1}^x f$ and $\sumr_{1}^x f$ is independent,
and if both are defined, their values can be different.

All results about right fractional sums in 
Sections~\ref{SecCharacterization} and
\ref{SecProducts} carry over to the case of left fractional sums.

\begin{lemma}[Left and Right Summation]
\label{LemLeftRight} \lineclear
We have
\[
\suml_{\nu=a}^b f(\nu) = \sumr_{\nu=-b}^{-a} f(-\nu)
\]
for all $a,b\in\C$ for which at least one of these sums exists.
\end{lemma}
\proof
Define a function $g(x):=f(-x)$. Suppose the right sum exists, so $g$ is a
right approximate polynomial and there are polynomials $q_n$ of fixed degree
$\sigma$ such that $|g(n+x)-q_n(n+x)|\to 0$ as $n\to+\infty$, for every $x$.
Then $f$ is a left approximate polynomial with approximating polynomials
$q_{-n}(-\nu)$.

By continued summation (Theorem~\ref{ThmLinearity}), then the definition, then
changing the sign of the summation index, then continued summation again, and
changing the sign of $n$, we obtain
\begin{eqnarray*}
&&
\sumr_{\nu=-b}^{-a} f(-\nu)
= -\sumr_{\nu=-a+1}^{-b-1} g(\nu)
\\
&=&
-\lim_{n\to\infty}\left(
\sumrl_{\nu=n-a+1}^{n-b-1}q_n(\nu)
+\sum_{\nu=1}^n (g(\nu-a) - g(\nu-b-1))\right)
\\
&=&
-\lim_{n\to\infty}\left(
\sumrl_{\nu=-n+b+1}^{-n+a-1} q_{n}(-\nu)
+\sum_{\nu=-n}^{-1} (g(-\nu-a) - g(-\nu-b-1))
\right)
\\
&=&
\lim_{n\to\infty}\left(
\sumrl_{\nu=-n+a}^{-n+b} q_{n}(-\nu)
+\sum_{\nu=0}^{-n-1} (g(-\nu-a) - g(-\nu-b-1))
\right)
\\
&=&
\lim_{n\to-\infty}\left(
\sumrl_{\nu=n+a}^{n+b}q_{-n}(-\nu)
+\sum_{\nu=1}^{n} (f(\nu+a-1) - f(\nu+b))\right)
= \suml_{\nu=a}^b f(\nu)
\end{eqnarray*}
as claimed.
\qed

\begin{theorem}[The Geometric Series]
\label{ThmGeometricSeries} \lineclear
For all $x\in \C$, we have
\begin{eqnarray*}
\sumr_{\nu=0}^x q^\nu =  \frac{q^{x+1}-1}{q-1} &&
\mbox{for $0\leq q<1$, \,\, and } \\
\suml_{\nu=0}^x q^\nu = \frac{q^{x+1}-1}{q-1} &&
\mbox{for $q>1$} \,\,.
\end{eqnarray*}
\end{theorem}
\proof
If $q\in[0,1)$, we get by resolving the definition
\[
\sumr_{\nu=0}^x q^\nu
= \sum_{\nu=1}^\infty \left(q^{\nu-1}-q^{\nu+x}\right)
= (1-q^{x+1}) \sum_{\nu=1}^\infty q^{\nu-1}
= \frac{q^{x+1}-1}{q-1}
\]
as claimed. The case $q>1$ is analogous.
\qed

\remark
The result carries over to the case of complex $q$ with $|q|<1$ 
resp.\ $|q|>1$, and
the proof is the same. Care has to be taken with branch cuts even in the sum
$\sumr_{\nu=0}^x q^\nu$ because of the occurrence of terms $q^{x+1+\nu}$. It
suffices to fix one branch of $q^{x+1}$ throughout the proof.

Now we show that the Binomial series makes sense even for
non-integer exponents. We use the general expression of binomial
coefficients
\[
{c \choose \nu}:=\frac{c!}{\nu! (c-\nu)!} \equiv
      \frac{\Gamma(c+1)}{\Gamma(\nu+1) \Gamma(c-\nu+1)} \,\, .
\]
Since $\Gamma$ is meromorphic in $\C$ without zeroes and with poles exactly at
non-positive integers, the binomial coefficient has a well-defined value in
$\Cbar$ for every $c,\nu\in\C$. We ignore the case that $c$ is a 
negative integer.
Then the binomial coefficient takes values in $\C$, and it vanishes 
exactly when
$\nu$ or $c-\nu$ are negative integers.

\begin{theorem}[The Binomial Series]
\label{ThmBinomialSeries} \lineclear
For all $c\in \C\setminus\{-1,-2,-3,\cdots\}$, we have
\begin{eqnarray*}
(1+x)^c=\sumr_{\nu=0}^c {c \choose \nu} x^\nu &&
\mbox{for all $x\in\C$ with $\vert x\vert <1$ \,\, and } \\
(1+x)^c=\suml_{\nu=0}^c {c \choose \nu} x^\nu &&
\mbox{for all $x\in\C$ with $\vert x\vert >1$} \,\,.
\end{eqnarray*}
\end{theorem}

\proof
For $\vert x \vert <1$, well-known estimates imply that the summand is right
summable with $\sigma=-\infty$. Therefore, we get by resolving the definition
\[
\sumr_{\nu=0}^c {c \choose \nu} x^\nu = {c \choose 0} +
     \sumr_{\nu=1}^c {c \choose \nu} x^\nu =
1+\sum_{\nu=1}^\infty \left ( {c \choose \nu} x^\nu -
      {c \choose {\nu+c}} x^{\nu+c} \right)
\,\,.
\]
The last binomial
coefficient always vanishes and we get
\[
\sumr_{\nu=0}^c {c \choose \nu} x^\nu =
     1+\sum_{\nu=1}^\infty {c \choose \nu} x^\nu = (1+x)^c
\]
as claimed. For $\vert x \vert >1$ we use Lemma~\ref{LemLeftRight}
and the first part to find
\[
\suml_{\nu=0}^c {c \choose \nu} x^\nu =
       \sumr_{\nu=-c}^0 {c \choose -\nu} x^{-\nu}
       = \sumr_{\nu=0}^c {c \choose c-\nu} x^{c-\nu}
       = x^c \sumr_{\nu=0}^c {c \choose \nu} \left(\frac 1 x\right)^{\nu}
= x^c \left(1+\frac 1 x\right)^c
\]
and the claim follows.
\qed

\section{Riemann and Hurwitz Zeta Functions}
\label{SecRiemannHurwitz}

The following Dirichlet series generalizes the Riemann $\zeta$
function and is known as the {\em Hurwitz $\zeta$ function}:
\begin{equation}
\zeta(s,x):=\sum_{\nu=0}^\infty \frac{1}{(\nu+x)^s} \,\,;
\label{EqHurwitzZeta}
\end{equation}
here $x$ is an arbitrary complex number but not a negative
integer or zero.
The series converges whenever $\Re(s)>1$. For the definition of
the powers $(\nu+x)^{-s}$, we use a branch cut at $\R^-$, and then
define the complex logarithm on $\R^-$ via continuity from above. This
way, the function $x\mapsto \zeta(s,x)$ is analytic on $\C\setminus(-\infty,0]$.

For every $x\in\C\sm(-\N)$, the Hurwitz $\zeta$ function extends to a meromorphic
function in $s$ with a single pole at $s=1$. One way to see this is via the
formula
\begin{eqnarray}
\frac{\partial}{\partial x} \zeta(s-1,x) = -(s-1)\, \zeta(s,x)
\quad\mbox{ for every }x\in\C\setminus(-\infty,0]
\label{EqHurwitzDiff}
\end{eqnarray}
from \cite[64:10:1]{SO} (which is easily verified directly).

The difference to the Riemann
$\zeta$ function is in the appearance of $x$, and in the fact that
summation starts with $\nu=0$; hence $\zeta(s,1)=\zeta(s)$.

For $x\in\R$, let $\lfloor x\rfloor$ denote the largest integer not
exceeding $x$.

\begin{lemma}[Polynomial Approximation of the Hurwitz Zeta Function]
\label{LemHurwitzProperty} \lineclear
For every compact $K\subset\C$ and every $s\in\C\sm\{1\}$, there is a
sequence of complex polynomials $p_n$ of degree $\sigma(s)$ such that
\begin{eqnarray*}
\left| \strut \zeta(s,n+z)-p_n(n+z) \right|
   &\longrightarrow& 0
\end{eqnarray*}
as $n\to\infty$, uniformly for $z\in K$, with
\[
\sigma(s) = \quad \left\{
\begin{array}{rl}
-\infty & \mbox{ if\/ $\Re(s)>1$, } \\
1+\lfloor\Re(-s)\rfloor & \mbox{ if\/ $\Re(s)\leq 1$ and $s\neq 1$.}
\end{array}
\right.
\]
In particular, $z\mapsto\zeta(s,z)$ is approximately polynomial of
degree $\sigma(s)$.
\end{lemma}

\proof
The case $\Re(s)>1$ is clear with $\sigma(s)=-\infty$, hence
$p_n(z)\equiv 0$. The case that $s\in(\C\sm\{1,0,-1,-2,\dots\})$
will be shown by induction on $s-N$ for $N\in\N$. We may suppose that
$K$ is a disk with center $z_0$. For every $s\neq 1$, the map $x\mapsto
\zeta(s,x)$ has no poles for $\Re(x)>0$; given $K$, we restrict to $n$
so all $z\in K$ satisfy $\Re(z+n)>0$.

For given $s$, let $p_n$ be approximating polynomials for $\zeta(s,.)$.
Let $P_n$ be polynomials with $P_n'=-(s-1)p_n$ and
$P_n(n+z_0)=\zeta(s-1,n+z_0)$; these are approximating polynomials
for $\zeta(s-1,.)$:
\begin{eqnarray*}
\zeta(s-1,n+z)-P_n(n+z)
&=&
\int_{z_0}^z \left( \strut
(\partial/\partial z')\zeta(s-1,n+z') - P'_n(n+z')\right)\,dz' \\
&=&
-(s-1)\int_{z_0}^z \left(\zeta(s,n+z')-p_n(n+z')\right) \, dz' \,\,.
\end{eqnarray*}

Finally, if $s\in\{0,-1,-2,\dots\}$, then $\zeta(s,z)$ is a polynomial
in $z$ \cite[64:4]{SO}.
\qed

\begin{corollary}[Extended Riemann-Hurwitz $\zeta$ Formula]
\label{CorRiemannHurwitz2}  \lineclear
For $x\in\C\setminus\{-1,-2,-3,...\}$ and $a\in\C\sm\{-1\}$, the
fractional sum $\sumr_{\nu=1}^x \nu^a$
exists and satisfies
\[
\sumr_{\nu=1}^x \nu^a = \zeta(-a)-\zeta(-a,x+1) \,\,.
\]
\end{corollary}
\proof
It is well known that
\begin{equation}
\zeta(s,x+1)=\zeta(s,x)-x^{-s}
\label{EqZetaDiff}
\end{equation}
for all $x\in\C\setminus\{-1,-2,-3,\cdots\}$ and $s\in\C\setminus\{1\}$. (For
$\Re(s)>1$, this follows directly from the definition. Since both sides are
holomorphic in
$s$ on $\C\setminus\{1\}$ for every $x$, the formula holds in general).

Let $S_a(x):=\zeta(-a)-\zeta(-a,x+1)$; then
$S_a(0)=\zeta(-a)-\zeta(-a,1)=0$ and
\[
     \Delta S_a(x)=S_a(x)-S_a(x-1)=x^a \,\,.
\]
Using  Lemma~\ref{LemHurwitzProperty}, the claim follows from
Lemma~\ref{LemSummationApprox}.
\qed

\remark
For $a=-1$, we have the formula
\begin{equation}
\sumr_{\nu=1}^x \frac 1 \nu = \sum_{\nu=1}^\infty \left(
        \frac 1 \nu - \frac 1 {\nu+x}
     \right)
= \gamma+ \frac{d}{dx}\ln\Gamma(x+1)=\gamma+\psi(x+1)
\label{EqHarmonicSeries}
\end{equation}
which converges whenever $x$ is not a negative integer. Here
$\gamma=0.577\dots$ is the Euler-Mascheroni constant and $\psi$ is
the so-called digamma function; see \cite[6.3.16.]{AbramowitzStegun}.
A special case is
\[
\sumr_{\nu=1}^{-\frac 1 2} \frac 1 \nu =
-2\ln 2
\]
which was noticed already by Euler \cite{Euler}.

In the following, it will be convenient to introduce the notation
\[
\up x a := \sumr_{\nu=1}^x \nu^a \,\,.
\]
By Corollary~\ref{CorRiemannHurwitz2} and its remark, this is defined for
$x\in\C\setminus\{-1,-2,-3,\cdots\}$ and $a\in\C$.
The value $x=+\infty$ is perfectly admissible when $\Re(a)<-1$ and yields
\[
\up\infty a=\zeta(-a)
\]
with the Riemann $\zeta$ function. For $x\in\N$, the function
$\up x a$ is known as the generalized harmonic series $H_x^{-a}$,
and $\up x a=1+2^a+3^a+\ldots+x^a$ are the first $x$ terms of the
Dirichlet series of $\zeta(-a)$ for $\Re(-a)>1$.

\begin{corollary}[Zeta Derivatives Identity]
\label{CorZetaDerivatives} \lineclear
For $x\in\C\setminus\{-1,-2,-3,...\}$, all $a\in\C\sm\{-1\}$ and
$b\in\N$, we have
\[
\sumr_{\nu=1}^x \nu^a\left(\ln\nu\right)^b =
     (-1)^b\left(\zeta^{(b)}(-a)-\zeta^{(b)}(-a,x+1)\right)
\]
where $\zeta^{(b)}=(\partial/\partial a)^b\zeta$.
\end{corollary}
\proof
This follows formally by differentiating Corollary~\ref{CorRiemannHurwitz2}
$b$ times with respect to $a$. We omit the proof that this formal
differentiation is allowed.
\qed

\begin{corollary}[Power Sums Up To $-1/2$]
\label{CorPowerSums-12} \lineclear
For all $n\in\N\setminus\{0\}$ and $a\in\C\setminus\{-1\}$, we have
\[
     \up {\left(-\frac 1 2\right)}a \equiv \sumr_{\nu=1}^{-\frac 1
2} \nu^a
     =\left(2-2^{-a}\right)\zeta(-a)
     \qquad \mbox{ and } \qquad
     \up {\left(-\frac 1 2\right)}{2n} \equiv \sumr_{\nu=1}^{-\frac
1 2} \nu^{2n}=0
     \,\, .
\]
\end{corollary}

\proof
This follows directly from Corollary~\ref{CorRiemannHurwitz2}
with the special values of the $\zeta$ functions at the
corresponding arguments.
\qed

\remark
For the second identity, there is also a more direct proof which
uses only properties of fractional sums instead those of the $\zeta$
functions. By Lemma~\ref{LemLeftRight} and Continued Summation,
we have for $z\neq 0$
\begin{equation}
\suml_{\nu=1}^{-1/2} \nu^z =
     \sumr_{\nu=1/2}^{-1} (-\nu)^z = -\sumr_{\nu=0}^{-1/2}(-\nu)^z
     = -(-1)^z \sumr_{\nu=1}^{-1/2} \nu^z - 0^z \,\,.
\label{EqSumOneHalf}
\end{equation}

For polynomials, left and right sum coincide trivially $\left(\strut \sumr p=\suml
p=\sumrl p\right)$, so if $z=2n$ for $n\in\N^+$, then
$\sumr_{\nu=1}^{-1/2}\nu^z=-\sumr_{\nu=1}^{-1/2}\nu^z$, which proves the claim.

This shows that if the left and right sum both exist and are
equal, there are interesting ways to manipulate those sums.
Moreover, we have the following corollary:

\begin{corollary}[Power Sums and Zeros of the Zeta Function]
\label{CorPowerSumsZeros} \lineclear
Modify the definition of $\nu^z$ so that $\nu^z=0$ for $\nu=0$. Then for
all $z\in\C\sm(2\Z+1)$ with $2^{-z}\neq 2$, we have
\[
     \sumr_{\nu=1}^{-1/2}\nu^z=\suml_{\nu=1}^{-1/2}\nu^z
     \quad\Longleftrightarrow\quad \zeta(-z)=0 \,\, .
\]
\end{corollary}

\remark
The definition $0^z:=0$ simply omits the $\nu=0$-term in $\suml$. In
the proof, this corresponds to the omission of the $(\nu=0)$-term in
(\ref{EqSumOneHalf}).

\proof
If $z\in\C$ is not an odd integer, then $-(-1)^z\neq 1$. If the two
sums exist and are equal, then they must vanish by
(\ref{EqSumOneHalf}). If in addition $2^{-z}\neq 2$, then $\zeta(-z)=0$
by Corollary~\ref{CorPowerSums-12}.
Conversely, if $\zeta(-z)=0$, then $\sumr_{\nu=1}^{-1/2}\nu^z=0$ by
Corollary~\ref{CorPowerSums-12} and $\suml_{\nu=1}^{-1/2}\nu^z=0$ by
(\ref{EqSumOneHalf}).
\qed

We are now going to show that series multiplication makes it possible to
evaluate multiple fractional sums of the powers. Due to Corollary~\ref{CorRiemannHurwitz2},
the following results can be considered as statements about the Hurwitz $\zeta$ function.

\begin{lemma}[Double Sums of Powers]
\label{LemQuadPowerGeneral} \lineclear
For arbitrary $a,b\in\C$ and $x\in\C\setminus\{-1,-2,-3,\cdots\}$,
we have
\begin{equation}
\sumr_{\nu=1}^x
\left( \up \nu a \cdot \nu^b +\up \nu b\cdot\nu^a \right)
=
(\up x a)(\up x b) + \up x{a+b}  \,\,.
\label{EqQuadPowerGeneral}
\end{equation}
\end{lemma}
\proof
The well-known asymptotic expansion
\[
   \zeta(s,x)=\frac{x^{1-s}}{s-1}+\frac{x^{-s}} 2 +\sum_{j=1}^{m-1}
   \frac{B_{2 j}\Gamma(2j+s-1)}{(2j)! \Gamma(s)} x^{-2j-s+1}
   +\mathcal{O}(x^{-2m-s-1})\,\,,
\]
mentioned for example in \cite{Adamchik}, shows by series
multiplication that the function $Z(x):=\zeta(-a,x+1)\cdot\zeta(-b,x+1)$ can
be approximated by a finite linear combination of monomials $\left(x^z\right)$ for $\Re(x)\to\infty$.
Since monomials are approximately polynomial (see Corollary~\ref{CorRiemannHurwitz2} together
with Lemma~\ref{LemSummationApprox}), so must be $Z$. By Corollary~\ref{CorRiemannHurwitz2}, it follows that the
function $x\mapsto\left(\strut\up x a\right)\cdot\left(\up x b\right)$ is an approximate polynomial.

Thus, $f(x)=x^a$ and $g(x)=x^b$ satisfy the conditions of Lemma~\ref{LemProdSums}. This yields
\begin{eqnarray*}
\left(\sumr_{\nu=1}^x \nu^a\right)\left(\sumr_{\nu=1}^x \nu^b\right)\equiv
(\up x a)(\up x b)
&=&
\sumr_{\nu=1}^x \left(
   -\nu^a \nu^b+\nu^a\sumr_{k=1}^\nu k^b + \nu^b \sumr_{k=1}^\nu k^a
\right)
\\
&=&
-\left(\up x {a+b}\right)+\sumr_{\nu=1}^x \left(
   \nu^a \cdot \up {\nu}b + \nu^b \cdot \up{\nu} a
\right)\,\,.
\end{eqnarray*}
\qed

\begin{corollary}[Double Power Sum]
\label{CorDoublePowerSum} \lineclear
For $a\in\C$ and $x\in\C\setminus\{-1,-2,-3,\cdots\}$, we have
\(\displaystyle
\sumr_{\nu=1}^x \up \nu a = \up x a (x+1) - \up x{a+1} \,\,.
\)
\end{corollary}

\proof
This follows from the second identity of Lemma~\ref{LemQuadPowerGeneral}
with
$b=0$, using $\nu^0=1$ and $\up\nu 0=\nu$:
\[
\sumr_{\nu=1}^x\up \nu a + \up x{a+1}
= \sumr_{\nu=1}^x \left(\up\nu a+\nu^{a+1}\right)
=\sumr_{\nu=1}^x\left(\up\nu a\cdot\nu^0+\up\nu 0\cdot\nu^a\right)
=(\up x a)(\up x 0)+\up x{a+0}
=\up x a(x+1) \,\,.
\]
\qed

We can use Corollary~\ref{CorDoublePowerSum} and
Lemma~\ref{LemQuadPowerGeneral}
iteratively with $b\in\N$ to compute arbitrary multiple sums of powers; a
lengthy calculation gives
\begin{eqnarray*}
\sumr_{k=1}^x k^a
&=&
\up x a
\\
\sumr_{l=1}^x \sumr_{k=1}^l k^a
&=&
\up x a (x+1) - \up x {a+1}
\\
\sumr_{m=1}^x\sumr_{l=1}^m\sumr_{k=1}^l k^a
&=&
\up x a \left( \frac{x^2}2+\frac 3 2 x +1
\right)-\up x {a+1} \left( x+\frac 3 2 \right) +\frac{\up x {a+2}} 2\,\,.
\hide{
}
\end{eqnarray*}
If we write ${\sumr}^{n+1}$ for the $n+1$-fold iterate of the summation operator, we
can write the general formula as
\[
   {\sumr}^{n+1} x^a=\frac 1 {n!} \sum_{\nu=0}^n
   \frac {\up x {a+\nu}} {\nu !} (-1)^\nu \frac{d^\nu}{dx^\nu}
   \prod_{k=1}^n (x+k)
\]
(this was observed by Michael Stoll). This is proved easily by induction: just show
that by application of the difference operator $\Delta$ on the right hand side, $n$
is replaced by $n-1$, and the claim follows by Lemma~\ref{LemSummationApprox}.
It is interesting to note that, after having defined $\up x
a$ for the first summation depending on $a$, the further iterated sums
can be expressed in terms of $\up x a$ only.

\section{Evaluation of Infinite Series, Products, and Limits}
\label{SecEvalInfinite}

Many infinite series can be evaluated quite intuitively
with the help of fractional sums. In this section, we give two examples.
First, consider the product
\begin{equation}
   P_1(x):=\lim_{n\to\infty}\prod_{k=1}^{2n+1}\left(1+\frac x k\right)^{k(-1)^{k+1}}
   \label{EqProdBorwein}
\end{equation}
which has been calculated by Borwein and Dykshoorn \cite{BD} and later by Adamchik \cite{Adamchik}.
We will now show that fractional sums allow to rederive a closed-form expression for $P_1$ in a straightforward way.

\begin{example}[Product by Borwein and Dykshoorn]
\label{ExBorweinDyk} \lineclear
For $P_1:\C\setminus\{-1,-2,-3,\ldots\}\to\C$ as defined in Equation~(\ref{EqProdBorwein}), we have
\[
   P_1(x)=2^{-\frac 1 {12}} \left(
      \frac{\Gamma\left(x+\frac 1 2\right)}{\Gamma(x+1)}
   \right)^{2 x} \exp\left({ -x-2\zeta'\left(-1,x+\frac 1 2\right)+2\zeta'\left(-1,x+1\right)-3\zeta'(-1)}\right)
   \,\,.
\]
\end{example}

\proof
The function $\nu\mapsto 2\nu\ln\left(1+\frac x \nu\right)$ is approximately polynomial
of degree $\sigma=0$, so by Definition~\ref{DefFracSummable}, using
$\lim_{n\to\infty} 2 n \ln\left(1+\frac x n\right)=2x$, we get
\begin{eqnarray*}
   &&\sumr_{\nu=1}^{-\frac 1 2}2\nu\ln\left(1+\frac x \nu\right)\\
   &=&\lim_{n\to\infty}\left[ 
       -\frac 1 2 \cdot 2x 
        +\sum_{k=1}^n 2 k \ln\left(1+\frac x k\right)
      -2\left(k-\frac 1 2\right)\ln\left(1+\frac x {k-\frac 1 2}\right)
   \right].
\end{eqnarray*}
Except for the summand $-x$, this equals the negative of the logarithm of $P_1$
as given in Equation~(\ref{EqProdBorwein}).
Now we use the following special case of Corollary~\ref{CorZetaDerivatives}
\begin{equation}
   \sumr_{\nu=1}^{-\frac 1 2}\nu^a\ln\nu=2^{-a}\zeta(-a)\ln 2-(2-2^{-a})\zeta'(-a)
   \qquad (a\in\C\setminus\{-1\})\,\,,
   \label{EqSumVlogV}
\end{equation}
together with index shifting, continued summation, Example~\ref{ExFactorial} and
Corollary~\ref{CorZetaDerivatives}:
\begin{eqnarray*}
   \ln P_1(x)&=& -x-\sumr_{\nu=1}^{-\frac 1 2} 2\nu\ln\left(1+\frac x \nu\right)=
   -x-2\sumr_{\nu=1}^{-\frac 1 2} \nu\ln\left(\frac{\nu+x}{\nu}\right)\\
   &=&-x-2\sumr_{\nu=1}^{-\frac 1 2}\nu\ln(\nu+x)+2\sumr_{\nu=1}^{-\frac 1 2}\nu\ln\nu\\
   &=&-x-2\sumr_{\nu=1+x}^{-\frac 1 2 +x} (\nu-x)\ln\nu -\frac {\ln 2}{12} -3 \zeta'(-1)\\
   &=&-x-2\sumr_{\nu=1+x}^{-\frac 1 2 +x} \nu\ln\nu +2 x \sumr_{\nu=1+x}^{-\frac 1 2 +x} \ln\nu
   -\frac {\ln 2}{12} -3 \zeta'(-1)\\
   &=&-x-2\left(\zeta'\left(-1,x+\frac 1 2\right)-\zeta'\left(-1,x+1\right)\right)\\
   && + 2 x \left(\ln\Gamma\left(x+\frac 1 2\right)-\ln\Gamma(x+1)\right)-\frac {\ln 2}{12}
   -3 \zeta'(-1)\,\,.
\end{eqnarray*}
The claim follows by exponentiation.
\qed

%
%


Now we give an example of a limit identity that seems to be new, as far as we know.
\begin{example}[Gamma Function Product]
\hide{
\[
   \lim_{n\to\infty}\left[
      e^{\frac n 8 (4n+1)} n^{-\frac 1 {16}-\frac 1 2 n(n+1)} (2 \pi)^{-\frac n 4}
      \prod_{k=1}^n \frac {\Gamma(k+1)^k}{\Gamma(k+\frac 1 2)^{(k-\frac 1 2)}}
   \right]
   =(2e)^{\frac 1 {24}} A^{\frac 3 4} \exp\left(-\frac{7\zeta(3)}{32\pi^2}\right)\,\,,
\]
}
\[  \lim_{n\to\infty}\left[     e^{\frac n 4(4n+1)}n^{-\frac 1 8 - n(n+1)}(2\pi)^{-\frac n 2}     \prod_{k=1}^{2n} \Gamma\left(1+\frac k 2\right)^{k(-1)^k}  \right]
=\sqrt[12]{2} \exp\left(\frac{5}{24} - \frac 3 2 \zeta'(-1) -\frac{7\zeta(3)}{16\pi^2}\right)
\,\,.
\]
\end{example}
\remark
Using the Glaisher-Kinkelin constant $A:=\exp\left(\frac  1 {12}-\zeta'(-1)\right)$, this identity
can also be written as
\[  \lim_{n\to\infty}\left[     e^{\frac n 4(4n+1)}n^{-\frac 1 8 - n(n+1)}(2\pi)^{-\frac n 2}     \prod_{k=1}^{2n} \Gamma\left(1+\frac k 2\right)^{k(-1)^k}  \right]
=(2 e)^{\frac 1 {12}} A^{\frac 3 2}\exp\left( -\frac{7\zeta(3)}{16\pi^2}\right)
\,\,.
\]
\proof
We use Lemma~\ref{LemProdSums} to multiply the fractional sums $\sumr_{\nu=1}^x \ln\nu$
and $\sumr_{\nu=1}^x v$:
\begin{eqnarray*}
   \left(\sumr_{\nu=1}^x \ln\nu\right)\cdot\left(\sumr_{\nu=1}^x \nu\right)
   =-\sumr_{\nu=1}^x \nu\ln\nu+\sumr_{\nu=1}^x \left(\ln\nu\sumr_{k=1}^\nu k\right)
   +\sumr_{\nu=1}^x \left(\nu \sumr_{k=1}^\nu \ln k\right)\,\,.
\end{eqnarray*}
Example~\ref{ExFactorial} yields
\[
   \ln(x!)\frac{x(x+1)} 2 = -\sumr_{\nu=1}^x \nu\ln\nu+\sumr_{\nu=1}^x \ln\nu
   \left(\frac {\nu^2} 2 + \frac\nu 2\right)+\sumr_{\nu=1}^x \nu\ln(\nu !)\,\,.
\]
Specializing $x=-\frac 1 2$ and using Corollary~\ref{CorZetaDerivatives}
and Equation~(\ref{EqSumVlogV}), we get
\begin{eqnarray*}
   \frac 1 2 \ln\pi \cdot\left(-\frac 1 8\right)
   &=&\frac 1 2 \sumr_{\nu=1}^{-\frac 1 2} \nu^2 \ln\nu-\frac 1 2 \sumr_{\nu=1}^{-\frac 1 2}\nu\ln\nu
   +\sumr_{\nu=1}^{-\frac 1 2}\nu\ln(\nu !)\\
   &=&\frac 1 2 \left(-\frac 7 4\right) \zeta'(-2)-\frac 1 2 \left(-\frac {\ln 2}{24}-\frac 3 2 \zeta'(-1)\right)
   +\sumr_{\nu=1}^{-\frac 1 2}\nu\ln(\nu !)\,\,.
\end{eqnarray*}
Rearranging gives
\[
   \sumr_{n=1}^{-\frac 1 2} n\ln(n!)=-\frac{\ln 2}{48}-\frac{\ln\pi}{16}-\frac 3 4 \zeta'(-1)
   +\frac 7 8 \zeta'(-2)\,\,.
\]
It is easily checked that $n\mapsto n\ln(n!)$ is approximately polynomial of degree $\sigma=2$.
Resolving the definition of this fractional sum and exponentiating, while using the
well-known identity $\zeta'(-2)=-\frac{\zeta(3)}{4\pi^2}$ and Stirling's approximation of $\ln(n!)$ for large $n$,
the claim follows.\qed

Many more identities can be produced in a similar spirit. Note that $-\frac 1 2$ is not the
only interesting fractional summation boundary; other examples, resulting in interesting
limit identities, include
\[
   \prodr_{n=\frac 1 4}^{-\frac 1 4} \left(n!\right)^n
   =\left(\frac{\Gamma\left(\frac 1 4\right)}{\Gamma\left(\frac 34\right)}\right)^{\frac 3 {32}}
   \exp\left(
      \zeta'\left(-2,\frac 1 4\right)-\frac{3\zeta(3)}{128\pi^2} -\frac G {4\pi}
   \right)\,\,,
\]
where $G=\sum_{k=0}^\infty \frac{(-1)^k}{(2k+1)^2}$ denotes Catalan's constant, or
\begin{eqnarray*}
   \prodr_{n=1}^{-\frac 1 2}\left(n!\right)^{\ln n}&=&\exp\left(\frac{\gamma^2}4+\frac{\gamma_1} 2
   -\frac{\pi^2}{48}+\frac{\ln^2 2}2-\frac{\ln^2\pi}8\right)\,\,,
\end{eqnarray*}
where $\gamma=0.577215\ldots$ and $\gamma_1=.072815\ldots$ are the
Euler-Mascheroni and Stieltjes constants.

\section{Questions}
\label{SecQuestions}

The paper ``On some strange summation formulas'' \cite{Strange} 
contains formulas
such as
\begin{equation}
\sum_{n=1}^\infty \frac{(-1)^n}{n^2}\cos\sqrt{b^2+\pi^2 n^2}
= \frac {\pi^2} 4 \left( \frac{\sin b} b -\frac {\cos b} 3
\right)
\label{Eq:Strange}
\end{equation}
and
\begin{equation}
\sum_{n=0}^\infty \frac {(-1)^n}{n+\frac 1
2}\frac{\sin\sqrt{b^2+\pi^2(n+1/2)^2}}{\sqrt{b^2+\pi^2(n+1/2)^2}}
=
\frac{\pi}{2}\frac{\sin b}{b}
\,\,.
\label{Eq:Strange2}
\end{equation}

The original proofs use Fourier transforms and properties of Bessel functions.
It is tempting to prove these (and many other) identities using 
fractional sums. In fact, there is a simple argument (see \cite[Section~9]{FracSumsProc}) 
to derive these identities  using termwise evaluation of fractional sums of power series:
\[
\sumr_{\nu=a}^b \left(\sum_{i=0}^\infty a_i \nu^i\right)
\stackrel{?}{=}
\sum_{i=0}^\infty \left(a_i \sumr_{\nu=a}^b \nu^i\right)
\,\,.
\]
We leave it as an open question to establish sufficient conditions when this is
justified, and thus to complete the proofs of (\ref{Eq:Strange}) and (\ref{Eq:Strange2}) by fractional sums.

A related open question is concerned with differentiation with respect to the
summation boundaries: under which conditions is the relation
\[
\frac{d}{dx}\sumr_{\nu=1}^x f(\nu)
\stackrel{?}{=}
c_f
+ \sumr_{\nu=1}^x f'(\nu)
\]
valid with a constant $c_f$? For example, we have (\ref{EqHarmonicSeries})
\[
\frac{d}{dx}\sumr_{\nu=1}^x \ln(\nu)
=
\frac{d}{dx}\ln\Gamma(x+1)
=-\gamma
+ \sumr_{\nu=1}^x \frac 1 \nu \,\,,
\]
where $\gamma\approx 0.577\dots$ is the Euler-Mascheroni constant.
This question arises for example when comparing fractional sums with the
Euler-Maclaurin summation formula.

\end{document}